\documentclass[brochure,english,12pt]{bourbaki}
\usepackage[matrix,arrow]{xy}
\usepackage{amssymb,amsfonts,amsmath,footnote,color,graphicx,url,tikz}
\usepackage[francais]{babel}
\usetikzlibrary{shapes,matrix,through,arrows,patterns}

\newcommand{\zg}{\mathring{g}}

\newcommand{\mcH}{{\mycal H}}

\newcommand{\bg}{\mathbf{g}}
\newcommand{\bR}{\mathbf{R}}
\newcommand{\bT}{\mathbf{T}}
\newcommand{\hyp}{{\mycal S}}
\newcommand{\mcM}{{\mycal M}}

\newcommand{\funnyr}{{r}}

\newcommand{\sourcesR}{R}

\newcommand{\hg}{{\hat g{}}}

\newcommand{\tbeta}{\tilde \beta}
\newcommand{\deltag}{h}

\newtheorem{Problem}[defi]{Problem}

\DeclareFontFamily{OT1}{rsfs}{}
\DeclareFontShape{OT1}{rsfs}{m}{n}{ <-7> rsfs5 <7-10> rsfs7 <10-> rsfs10}{}
\DeclareMathAlphabet{\mycal}{OT1}{rsfs}{m}{n}

{\catcode `\@=11 \global\let\AddToReset=\@addtoreset}
\AddToReset{equation}{section}
\renewcommand{\theequation}{\thesection.\arabic{equation}}
\newcounter{mnotecount}[section]

\renewcommand{\themnotecount}{\thesection.\arabic{mnotecount}}

\newcommand{\mnote}[1]
{\protect{\stepcounter{mnotecount}}$^{\mbox{\footnotesize
$
\bullet$\themnotecount}}$ \marginpar{
\raggedright\tiny\em
$\!\!\!\!\!\!\,\bullet$\themnotecount: #1} }

\newcommand{\const}{\textrm{const.}}

 \newcommand{\bel}[1]{\begin{equation}\label{#1}}
 \newcommand{\bea}{\begin{eqnarray}}
 \newcommand{\bean}{\begin{eqnarray}\nonumber}
 \newcommand{\beal}[1]{\begin{eqnarray}\label{#1}}
\newcommand{\eeal}[1]{\label{#1}\end{eqnarray}}
\newcommand{\eea}{\end{eqnarray}}
\newcommand{\be}{\begin{equation}}
\newcommand{\eeq}{\end{equation}}
\newcommand{\ee}{\end{equation}}
\newcommand{\beqa}{\begin{eqnarray}}
\newcommand{\eeqa}{\end{eqnarray}}
\newcommand{\beqan}{\begin{eqnarray*}}
\newcommand{\eeqan}{\end{eqnarray*}}
\newcommand{\ba}{\begin{array}}
\newcommand{\ea}{\end{array}}

\newcommand{\eq}[1]{(\ref{#1})}

\newcommand{\trg}{{\mathrm{tr}_g}}

\newcommand{\R}{\mathbb R}

\newcommand{\N}{\mathbb N}

\date{Novembre 2016}
\bbkannee{69\`eme ann\'ee, 2016-2017}
\bbknumero{1120\hfill Preprint UWThPh 2016-18}
\title{Anti-gravity \`a la Carlotto-Schoen}
\subtitle{after Carlotto and Schoen}
\author{Piotr T. CHRU\'SCIEL}
\address{University of Vienna\\
Erwin Schr\"odinger Institute and Faculty of Physics\\
Boltzmanngasse 5\\
A-1090 Vienna, Austria}
\email{Piotr.Chrusciel@univie.ac.at}

\begin{document}
\maketitle

\noindent{\bf INTRODUCTION}

In~\cite{CarlottoSchoen} Carlotto and Schoen show that gravitational fields can be used to shield gravitational fields. {That is to say, one can produce spacetime regions extending to infinity where no gravitational forces are felt whatsoever, by manipulating the gravitational field around these regions.}  A sound-bite version of the result reads:

\medskip

\noindent{{\sc Theorem A} (Carlotto \& Schoen~~\cite{CarlottoSchoen}). ---
    \emph{Given an asymptotically flat  initial data set for vacuum Einstein equations there exist cones and asymptotically flat vacuum initial data which coincide with the original ones inside the cones and are   Minkowskian outside slightly larger cones, see Figure~\ref{F11IX16.1}.}
%
%
\vspace{-.3cm}
\begin{figure}[h]
  \centering
    \includegraphics[scale=.2]{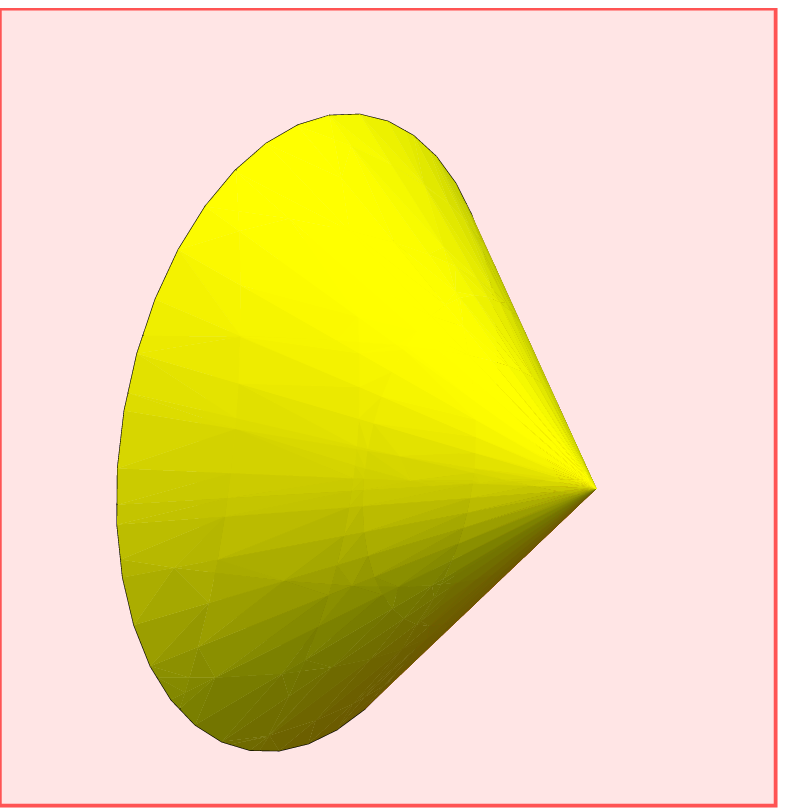}
    \includegraphics[scale=.25]{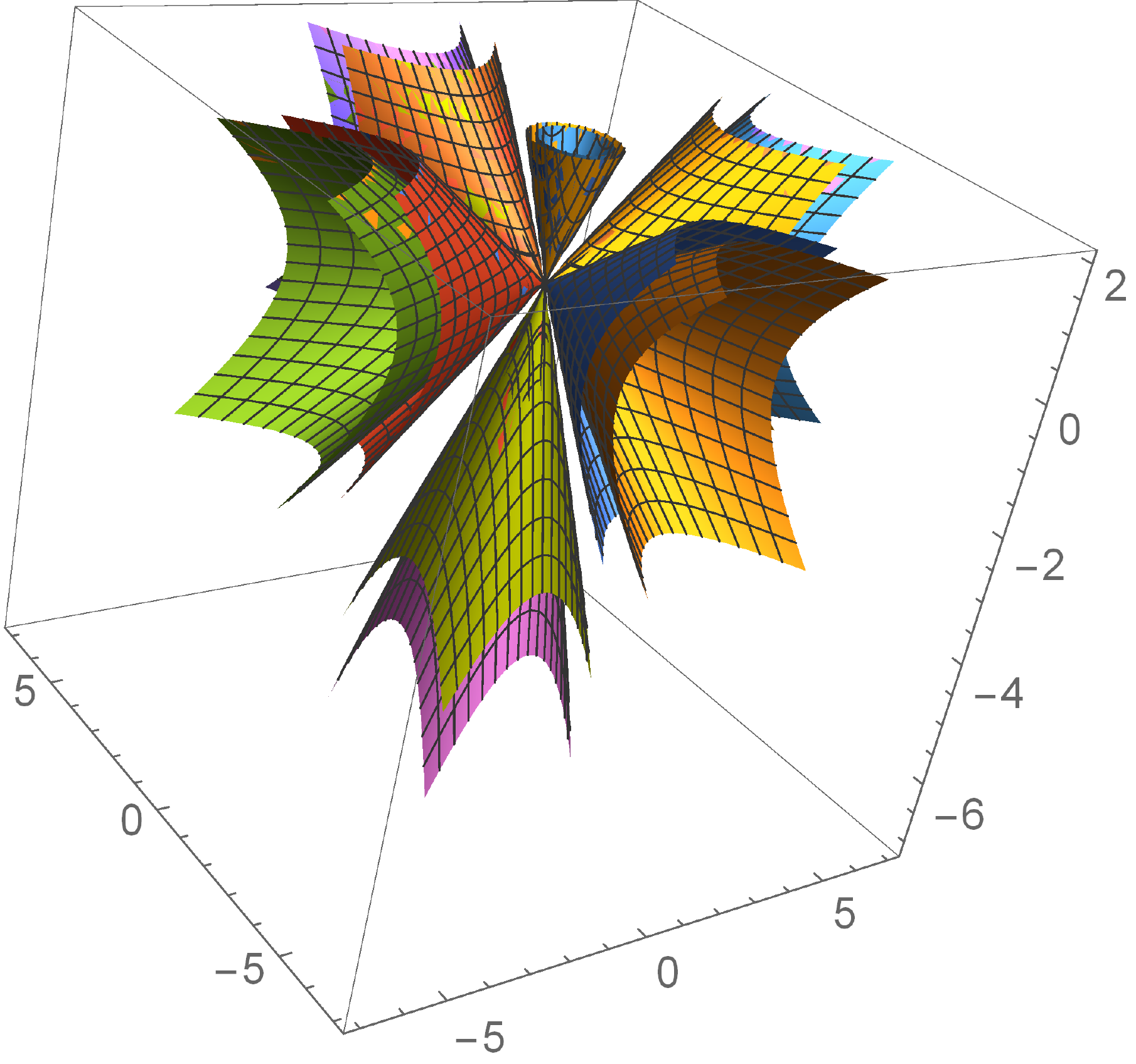}
    \vspace{-.3cm}
\caption{Left picture: {The new initial data are Minkowskian outside the larger cone, and coincide with the original ones inside the smaller one.} The construction can also be carried out the other way round, with Minkowskian data inside the smaller cone and the original ones outside the larger cone. Both cones extend to infinity, and their tips are located very far in the asymptotically flat region.
Right picture: Iterating the construction, one can embed any finite number of distinct initial data sets into Minkowskian data, or paste-in Minkowskian data inside several cones into a given data set.
    \label{F11IX16.1}}
\end{figure}

 {Actually, the result is true for all cones with preassigned axis and pair of apertures   provided the vertex is shifted sufficiently far away in the asymptotically flat regions.}

In the associated spacetimes $(\mcM,\bg)$ the  metric coincides with the Minkowski metric within the domain of dependence, which we will denote by $\mycal D$, of the complement of the larger cones, which forms an open subset of $\mcM$;  {we return to this in Section~\ref{s12X16.2} below}. Physical objects in $\mycal D$ do not feel any gravitational fields. The Carlotto-Schoen gluing has effectively switched-off any gravitational effects in this region. This has been achieved  by manipulating vacuum gravitational fields only.

In this \emph{s\'eminaire} I will describe the context, define the notions used, highlight the key elements of the  proof of Theorem~A, and discuss some  developments.

\section{The context}

\subsection{Newtonian gravity}

Newtonian gravity is  concerned with a gravitational potential field $\phi$, which solves the equation
\bel{11IX16.1}
 \Delta \phi = 4 \pi G \rho
 \,,
\ee
where $\Delta$ is the Laplace operator in an Euclidean $\R^3 $ and $G$ is Newton's constant. Up to conventions on signs, proportionality factors, and units, $\rho$ is the matter density, which is not allowed to be negative. Isolated systems are defined by the requirement that both $\rho$ and $\phi$ decay to zero as one recedes to infinity.

Freely falling bodies experience an acceleration proportional to the gradient of $\phi$. So no gravitational forces exist in those regions where $
\phi$ is constant.

Suppose that $\rho$ has support contained in a compact set $K$, and that $\phi$ is constant on an open set $\Omega$.  Since solutions of
\eq{11IX16.1} are analytic on $\R^3\setminus K$, $\phi$ is constant on any connected component of  $\R^3\setminus K$ which meets
$\Omega$. We conclude that if $\Omega$ extends to infinity, then $\phi $ vanishes at all large distances. This implies, for all sufficiently large spheres $S(R)$,
$$
 0 =
 \int_{S(R)} \nabla \phi \cdot n  \, d^2S =
 \int_{B(R)} \Delta\phi \, d^3V  =
  4\pi G
 \int_{B(R)} \rho \, d^3V
  \,.
$$
Since $\rho$ is non-negative, we conclude that $\rho \equiv 0$. Equivalently, for isolated systems with compact sources, \emph{Newtonian gravity cannot be screened away on open sets extending to large distances.}

The striking discovery of Carlotto and Schoen is, that this can be done in Einsteinian gravity.

The Newtonian argument above fails if matter with negative density is allowed. It should therefore be emphasised that the Carlotto-Schoen construction is done by manipulating vacuum initial data, without involvement of matter fields.

\subsection{Einsteinian gravity, general relativistic initial data sets}

\emph{Mathematical general relativity} is born around 1952  with the breakthrough paper of Yvonne Choquet-Bruhat~\cite{ChBActa}, showing that Einstein's field equations,
\bel{11IX16.1+}
 \bR_{\mu\nu} - \frac 12 \bR \bg_{\mu\nu}+
 \Lambda \bg_{\mu\nu} = \frac{8 \pi G}{c^4} \bT_{\mu\nu}
 \,,
\ee
admit a well posed Cauchy problem. Here $\bR_{\mu\nu}$ is the Ricci tensor of the spacetime metric $\bg$, $\bR$ its Ricci scalar,  $\bT_{\mu\nu}$ the energy-momentum tensor of matter fields, $\Lambda$   the ``cosmological constant'', $G$ is Newton's constant  as before, and $c$ is the speed of light. In vacuum $\bT_{\mu\nu}$ vanishes, in which case
\eq{11IX16.1+} is equivalent to the requirement that $\bg$ be Ricci-flat when moreover the vanishing of $\Lambda$ is imposed. The notation $\bg_{\mu\nu}$ indicates that the metric   $\bg$ is a two-covariant tensor field, similarly for $\bT_{\mu\nu}$, etc.

The geometric initial data for the four-dimensional vacuum Einstein equations are a triple $(\hyp,g,K)$, where $(\hyp,g)$ is a three-dimensional Riemannian manifold and $K$ is a symmetric two-covariant tensor field on $\hyp$. One should think of $\hyp$ as a space-like hypersurface in the vacuum Lorentzian spacetime $(\mcM,\bg)$, then $g$ is the metric induced by $
\bg$ on $\hyp$, and $K$ is the second fundamental form (``extrinsic curvature tensor'') of $\hyp$ in  $(\mcM,\bg)$.

 It has already been recognized in 1927 by Darmois~\cite{Darmois} that $(g,K)$ are not arbitrarily specifiable, but have to satisfy a set of constraint equations,
\beal{11IX16.3}
&
 R = |K|^2 - (\trg K)^2 + 2 \mu + 2 \Lambda
 \,,
 &
\\
&
 D ^i (K_{ij}-  \trg K g_{ij} ) =   J_j
 \,,
 &
\eeal{11IX16.4}
where $\mu = \frac {8 \pi G}{c^4} \bT_{\mu\nu} n^\mu n^\nu$ is the matter energy density on $\hyp$ and $J_j = \frac {8 \pi G}{c^4} \bT_{\mu j} n^\mu $ the matter momentum vector, with $n^\mu$ being the unit normal to $\hyp$ in $(\mcM,\bg)$. The requirement of positivity of energy of physical matter fields translates into the \emph{dominant energy condition}:
\bel{11IX16.5}
 \mu \ge |J|_g
 \,,
\ee
where of course $|J|_g \equiv \sqrt{g(J,J)} \equiv \sqrt{g_{ij} J^i J^j}$ (summation convention). In particular $\mu$ should be non-negative.

The constraint equations are the source of many headaches in mathematical and numerical general relativity. On the other hand, together with the energy condition, they are the   source of beautiful mathematical results%
\footnote{\label{Fn12X16.1}The reader is invited to consult~\cite{BartnikIsenberg,YCB:GRbook,CGP,CorvinoPollack} for more details and further references.}
such as the positive energy theorems~\cite{SchoenYau81,Witten81}, the Penrose inequality~\cite{Bray:preparation2,huisken:ilmanen:penrose}, the Corvino-Schoen~\cite{Corvino,CorvinoSchoen2} or the Carlotto-Schoen gluings.

\subsection{Asymptotic flatness}

Initial data for general relativistic isolated systems are typically modelled by \emph{asymptotically flat}   data with vanishing cosmological constant $\Lambda$. Actually, astrophysical observations indicate that $\Lambda$ is positive. However, for the purpose of observing nearby stars, or for our stellar system, the corrections arising from $\Lambda$ are negligible, they only become important at cosmological scales.

The class of asymptotically flat systems should obviously include the Schwarzschild black holes. In those, on  the usual slicing by $t=\const$ hypersurfaces it holds that $K_{ij}\equiv 0$ and
\bel{12IX16.1}
 g_{ij} = \left(1 + \frac m {2r}\right)^4 \delta_{ij} + O (r^{-2})
 \,,
\ee
 in spacetime dimension four, or
\bel{8X16.1}
 g_{ij} = \left(1 + \frac m {2r^{n-2}}\right)^\frac4{ n-2} \delta_{ij} + O (r^{-(n-1)})
 \,,
\ee
in general spacetime dimension $n+1$. Here $\delta_{ij}$ denotes the Euclidean metric in manifestly flat coordinates. The asymptotics \eq{12IX16.1} is often referred to as \emph{Schwarzschildean}, and the parameter $m$ is called the ADM mass of the metric. Now, one can obtain initial data with non-vanishing total momentum by taking Lorentz-transformed slices in the Schwarzschild spacetime. This leads to initial data sets satisfying
\beal{12IX16.2}
 &
\partial_{i_1}\cdots \partial_{i_{\ell}}( g_{ij}- \delta_{ij})=O (r^{-\alpha-\ell})
\,,
&
\\
 &
\partial_{i_1}\cdots \partial_{i_k} K _{ij} =O (r^{-\alpha-k-1})
 \,,
 &
\eeal{+12IX16.2}
with $\alpha=n-2$, for any $k,\ell\in \N$. Metrics $g$ satisfying \eq{12IX16.2} will be called \emph{asymptotically Euclidean}.

The flexibility of choosing $\alpha \in (0,n-2)$ in the definition of asymptotic flatness \eq{12IX16.2}-\eq{+12IX16.2}, as well as $k,\ell$ smaller than some threshold, is necessary in Theorem~A. Indeed, the new initial data constructed there are \emph{not expected} to satisfy
\eq{12IX16.2} with $\alpha=n-2$. It would be of interest to settle the question, whether or not this is really case.

There does not appear to be any justification for the Schwarzschildean threshold $\alpha=n-2$ other than historical. On the other hand, the threshold
\bel{8X16.3}
 \alpha = (n-2)/2
\ee
appears naturally as the optimal threshold for a well-defined total energy-momentum of the initial data set. This has been first discussed in~\cite{Chremark,ChErice,Chmass,DenisovSolovev}, compare~\cite{Bartnik}.

\subsection{Time-symmetric initial data and the Riemannian context}

Initial data are called \emph{time-symmetric} when $K_{ij}\equiv 0$. In this case, and assuming vacuum, the \emph{vector constraint equation}  \eq{11IX16.4} is trivially satisfied, while the \emph{scalar constraint equation} \eq{11IX16.3} becomes the requirement that $(\hyp,g)$ has constant scalar curvature $R$:
\beal{9X16.1}
&
 R =   2 \Lambda
 \,.
 &
\eeal{11IX16.4+}
In particular $(\hyp,g)$ should be scalar-flat when $\Lambda=0$.
(The ``time-symmetric'' terminology reflects the fact that a suitable reflection across $\hyp$ in the
associated spacetime is an isometry.)
So all statements about vacuum initial data translate immediately into statements concerning scalar-flat Riemannian manifolds. For example, the following statement is a special case of Theorem~A :

\begin{theo}[Carlotto \& Schoen]
  \label{T11IX16.2} Given a scalar-flat asymptotically Euclidean metric $g$ there exist cones and scalar-flat asymptotically Euclidean metrics which coincide with $g$ {inside of the cones and are flat outside slightly larger cones.}
\end{theo}

This theorem was certainly one of the motivations for the proof of Theorem~A. Indeed, the question
of existence of non-trivial, scalar-flat, asymptotically flat metrics $\hat g$ which are exactly flat in a half-space arises when studying complete, non-compact minimal hypersurfaces. Indeed, if $\hat g$ is such a metric,
then all hyperplanes lying in the flat half-space minimize area under compactly supported deformations which do not extend into the non-flat region. So Theorem~\ref{T11IX16.2} shows that such metrics $\hat g$ actually exist. This should be contrasted with the following beautiful result of Chodosh and Eichmair~\cite{CarlottoChodoshEichmair}, which shows that \emph{minimality under all compactly supported perturbations} implies flatness:

\begin{theo}[Chodosh, Eichmair]
  \label{T11IX16.3}
  The only asymptotically Euclidean three-dimensional manifold with non-negative scalar curvature that contains a complete non-compact embedded  surface $S$ which is a (component of the) boundary of some properly embedded full-dimensional submanifold of $(M, g)$ and is
area-minimizing under compactly supported deformations is flat $\R^3$, and $S$ is a flat plane.
\end{theo}

{The above was preceded by a related rigidity result of Carlotto~\cite{Carlotto}:}

\begin{theo}[Carlotto]
  \label{T11IX16.4}
  {Let $(M,g)$ be a complete, three-dimensional, asymptotically Schwarzschildean Riemannian manifold with non-negative scalar curvature. If $M$ contains a complete, properly embedded, stable minimal surface $S$, then $(M,g)$ is the Euclidean space and $S$ is a flat plane.}
\end{theo}

Such results immediately imply non-compactness for sequences of solutions of the Plateau problem with a diverging sequence of boundaries. We note that compactness results in this spirit play a key role in the Schoen \& Yau proof of the positive energy theorem. One could likewise imagine that convergence of such sequences of solutions of the Plateau problem could provide a tool to study stationary black hole solutions, but no arguments in such a spirit have been successfully implemented so far.

\subsection{Localised scalar curvature}
 \label{s12X16.1}

It is an immediate consequence of the positive energy theorem that,  for complete asymptotically Euclidean manifolds with non-negative scalar curvature,
$$
 \mbox{\emph{curvature cannot be localised in a compact set}.}
$$
In other words, a flat region cannot enclose a non-flat one. A similar statement applies for general relativistic initial data sets satisfying the dominant energy condition~\eq{11IX16.5}. Indeed a  metric which is flat outside of a compact set would have zero total mass and hence would be flat everywhere by the rigidity-part of the positive energy theorem~\cite{Bartnik,SchoenYau79b}. Similarly for initial data sets~\cite{ChBeig1,ChMaerten}.

One would then like to know  \emph{how much flatness can a non-trivial initial data set carry?} This question provided another motivation for Theorem~A, which shows that non-flatness can be localised within cones.

A previous family of non-trivial asymptotically Euclidean scalar-flat metrics containing flat regions is provided by the \emph{quasi-spherical metrics} of Bartnik~\cite{Bartnik93}. In Bartnik's examples flatness can be localised within balls.

In~\cite{CarlottoSchoen} it is noticed that \emph{non-flat} regions \emph{cannot} be sandwiched between parallel planes. This follows immediately from the following formula for ADM mass due to Beig~\cite{BeigKomar} (compare~\cite{ashtekar:hansen,Chremark,Herzlich:mass}),
\bel{12X16.1}
 m = \lim_{R\to\infty}\frac 1 { 16 \pi} \int_{r=R} G_{ij} x^i n^j dS
 \,,
\ee
where $G_{ij}:= R_{ij} - \frac 12 R g_{ij}$ is the Einstein tensor, and $x^i$ is the coordinate vector in the asymptotically Euclidean coordinate system. Indeed, \eq{12X16.1} together with a non-zero mass and asymptotic flatness imply that the region where the Ricci tensor has to have non-trivial angular extent as one recedes to infinity. But this is not the case for a region sandwiched between two parallel planes.

\subsection{Focussed gravitational waves?}
 \label{s12X16.2}

Consider non-trivial asymptotically flat Carlotto-Schoen initial data set, at $t=0$, with the ``non-Minkowskianity'' localised in a cone with vertex at $\vec a$, axis $\vec i$  and aperture $\theta>0$, which we will denote by  $C({\vec a,\vec i, \theta})$.
It follows from~\cite{BieriJDG} that the associated vacuum space-time will exist globally when the data are small enough, in a norm compatible with the Carlotto-Schoen setting.

One can think of the associated vacuum space-time as describing a gravitational wave localised, at $t=0$,   in an angular sector of opening angle $2 \theta$ and direction defined by the vector $\vec i$.
The discussion of the previous section shows that $\theta$ can be made as small as desired but cannot be zero, so that all  such solutions must have non-trivial angular extent.

It is of interest to enquire how Carlotto-Schoen solutions evolve in time. In what follows we assume that $\theta<\pi/2$. Standard results on Einstein equations show that the boundary enclosing the non-trivial region travels outwards no faster than the speed of light $c$.  This, together with elementary geometry shows that at time $t$ the space-time metric will certainly be flat outside of a cone
\bel{15X16.1}
 C\Big({\vec a - \frac {c } {\sin(\theta)} t \,\vec i, \vec i, \theta}\Big)
 \,.
\ee

The reader will note that the tip of the cone \eq{15X16.1} travels faster than light, which is an artefact of the rough estimate. A more careful inspection near the tip of the cone shows that the domain of dependence at time $t$ consists of a cone of aperture $\theta$ spanned tangentally on the boundary of a sphere of radius $t$ as shown in Figure~\ref{F15X16.1}.
\begin{figure}[th]
  \centering
    \includegraphics[scale=.25]{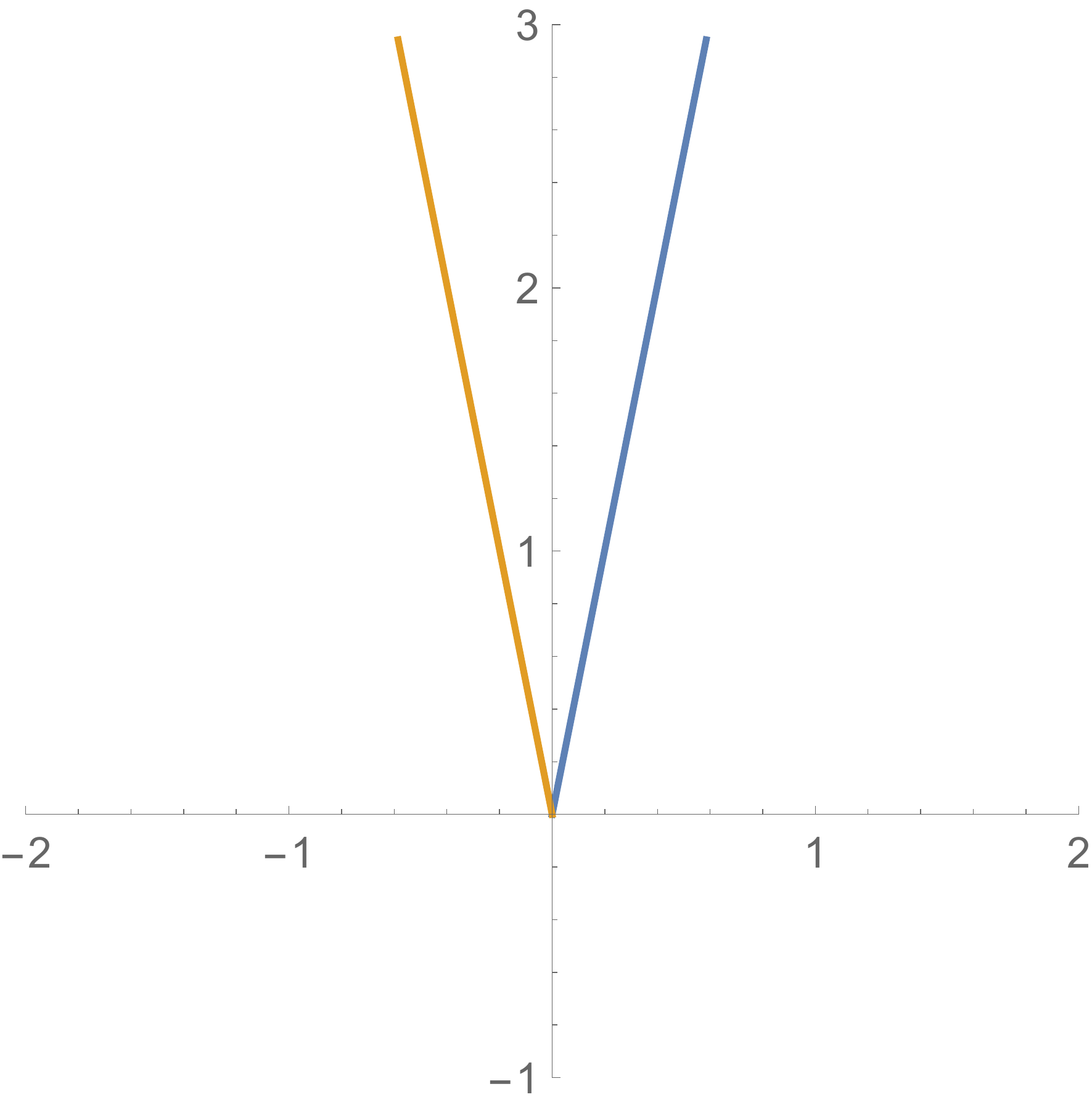}$t=0\longrightarrow t=1$
    \includegraphics[scale=.25]{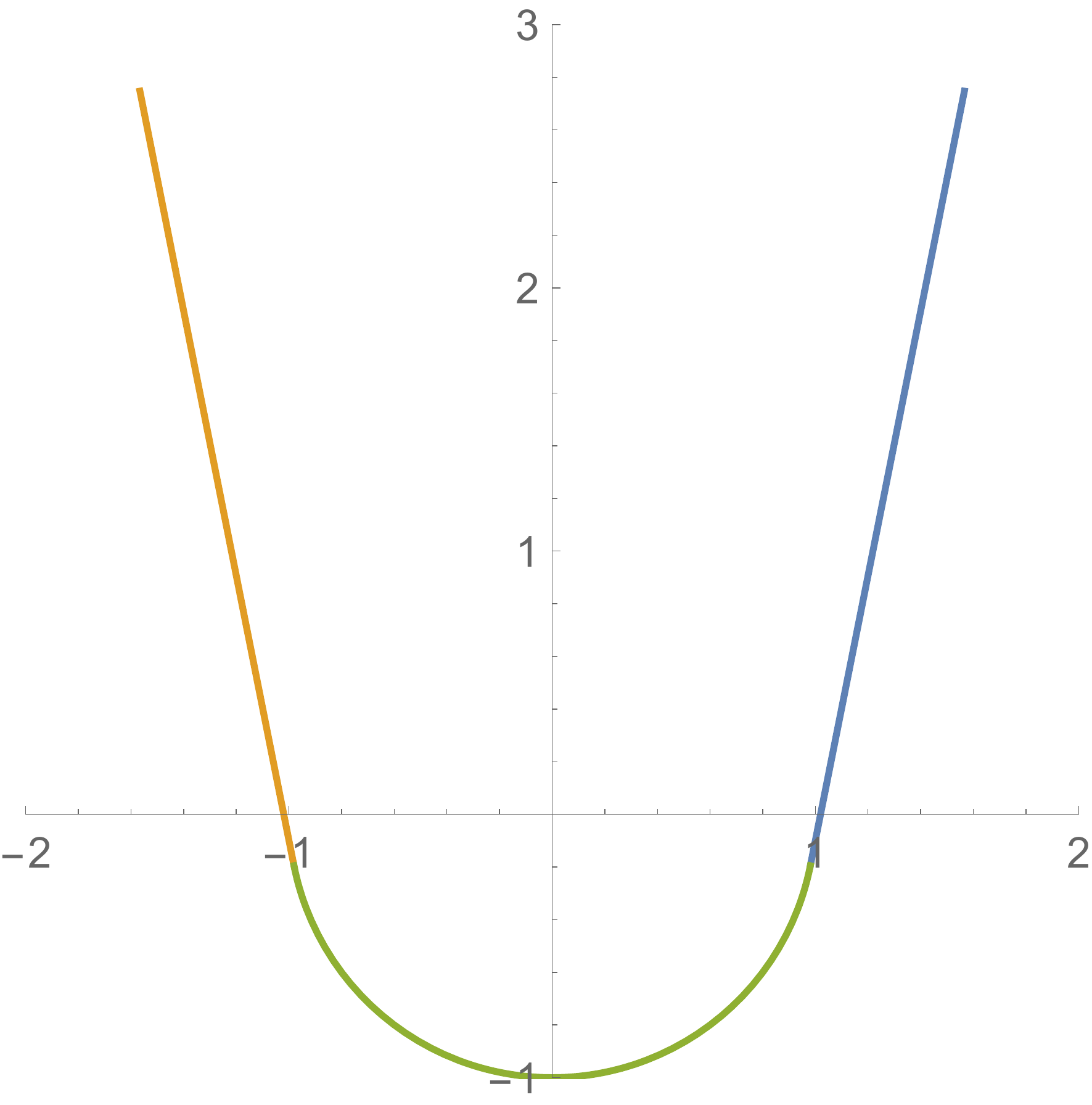}
\caption{Left graph:  The cone in the figure represents the exterior boundary of $\Omega$ at $t=0$, with initial data Minkowskian below the graph.
Right graph:  At $t=1$, the evolved space-time metric is Minkowskian below the graph. The three-dimensional picture is obtained by rotating the graphs around the vertical axis.
    \label{F15X16.1}}
\end{figure}
In any case the angular opening of the wave remains constant on slices of constant time. However, the wave is likely to spread and meet  all generators of null infinity, but not before  the intersection of $
\mycal S$ with the light-cone of the origin $(t=0,
\vec x=0)$ is reached.

\subsection{Many-body problem}
 \label{ss16X16.11}

Given two  initial data sets $(\Omega_a,g_a,K_a)$, $a=1,2$, the question arises, whether one can  find a new initial data set which contains both? An answer to this is not known in full generality. However, the Carlotto-Schoen construction gives  a positive answer to this question when the original initial data  are part of asymptotically flat initial data $(\hyp_a,g_a,K_a)$, provided that the sets $\Omega_a\subset \hyp_a$ can be enclosed in cones which do not intersect after ``small angular fattenings''. Indeed, one can then apply the deformation of Theorem~A to each original data set to new initial data $(\hyp_a,\hat g_a,\hat K_a)$ which coincide with the original ones on $\Omega_a$ and are Minkowskian outside the fattened cones. But then one can superpose the resulting initial data sets in the Minkowskian region, as shown in Figure~\ref{F16X16.1}.
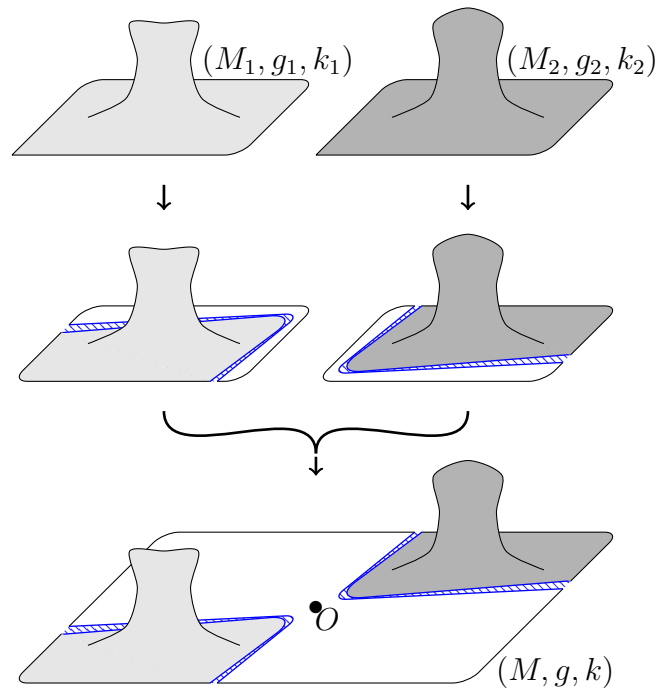
\begin{figure}[ht!]
\begin{center}
	\begin{tikzpicture}
	\begin{scope}[scale=0.5]
	
	\begin{scope}[shift={(0, 6)}]
	\begin{scope}[shift={(-4, 0)}]
	\fill [black!10, draw=black, rounded corners=2mm] plot coordinates {(-2, -1) (0, 1) (6, 1) (4, -1) (-2, -1)};
	\fill [black!10, draw=black] plot [smooth] coordinates {(0, 0) (1, 0.5) (1.3, 1.5) (1.1, 2.5) (2, 2.5) (2.9, 2.5) (2.7, 1.5) (3, 0.5) (4, 0)};
	\node at (5, 1.5) {$(M_1, g_1, k_1)$};
	\end{scope}
	
	\begin{scope}[shift={(4, 0)}]
	\fill [black!30, draw=black, rounded corners=2mm] plot coordinates {(-2, -1) (0, 1) (6, 1) (4, -1) (-2, -1)};
	\fill [black!30, draw=black] plot [smooth] coordinates {(0, 0) (1, 0.5) (1.25, 1.5) (1.1, 2.5) (2, 2.9) (2.9, 2.5) (2.75, 1.5) (3, 0.5) (4, 0)};
	\node at (5, 1.5) {$(M_2, g_2, k_2)$};
	\end{scope}
	
	\node at (-2, -1.5) (M1) {};
	\node at (6, -1.5) (M2) {};
	\node at (-2, -2.75) (M1tempup) {};
	\node at (6, -2.75) (M2tempup) {};
	\path[every node/.style={font=\sffamily\small}, ->, line width=1pt] (M1) edge [out=270, in=90] (M1tempup);
	\path[every node/.style={font=\sffamily\small}, ->, line width=1pt] (M2) edge [out=270, in=90] (M2tempup);
	
	\end{scope}
	
	\begin{scope}[shift={(-4, 0)}]
	\fill [pattern color=blue!80, pattern=north west lines, draw=blue, line width=0.5pt, rounded corners=2mm] plot coordinates {(-0.5, 0.5) (5.6, 0.8) (3.4, -1)};
	\fill [pattern color=blue!80, pattern=north west lines, draw=none, rounded corners=2mm] plot coordinates {(3.4, -1) (-2, -1) (-0.5, 0.5)};
	\fill [black!10, draw=blue, line width=0.5pt, rounded corners=4mm] plot coordinates {(-0.7, 0.29) (5.6, 0.8) (3.2, -1.01)};
	\fill [black!10, draw=black, rounded corners=2mm] plot coordinates {(3.2, -1) (-2, -1) (-0.7, 0.3)};
	\draw [black, rounded corners=2mm] plot coordinates {(-0.5, 0.5) (0, 1) (6, 1) (4, -1) (3.4, -1)};
	\fill [black!10, draw=black] plot [smooth] coordinates {(0, 0) (1, 0.5) (1.3, 1.5) (1.1, 2.5) (2, 2.5) (2.9, 2.5) (2.7, 1.5) (3, 0.5) (4, 0)};
	\end{scope}
	
	\begin{scope}[shift={(4, 0)}]
	\fill [pattern color=blue!80, pattern=north west lines, draw=blue, line width=0.5pt, rounded corners=2mm] plot coordinates {(0.6, 1) (-1.6, -0.8) (4.5, -0.5)};
	\fill [pattern color=blue!80, pattern=north west lines, draw=none, rounded corners=2mm] plot coordinates {(4.5, -0.5) (6, 1) (0.6, 1)};
	\fill [black!30, draw=blue, line width=0.5pt, rounded corners=4mm] plot coordinates {(0.8, 1.01) (-1.6, -0.8) (4.7, -0.29)};
	\fill [black!30, draw=black, rounded corners=2mm] plot coordinates {(4.7, -0.3) (6, 1) (0.8, 1)};
	\draw [black, rounded corners=2mm] plot coordinates {(0.6, 1) (0, 1) (-2, -1) (4, -1) (4.5, -0.5)};
	\fill [black!30, draw=black] plot [smooth] coordinates {(0, 0) (1, 0.5) (1.25, 1.5) (1.1, 2.5) (2, 2.9) (2.9, 2.5) (2.75, 1.5) (3, 0.5) (4, 0)};
	\end{scope}
	
	\node at (-2, -1.5) (M1temp) {};
	\node at (6, -1.5) (M2temp) {};
	\node at (2, -3.2) (M) {};
	\path[every node/.style={font=\sffamily\small}, line width=1pt] (M1temp) edge [out=270, in=90] (M);
	\path[every node/.style={font=\sffamily\small}, line width=1pt] (M2temp) edge [out=270, in=90] (M);
	\draw [black, ->, line width=1pt] (2, -3) -- (2, -3.5);
	
	\begin{scope}[shift={(0, -7)}]
	\begin{scope}[shift={(-4, -1)}]
	\fill [pattern color=blue!80, pattern=north west lines, draw=blue, line width=0.5pt, rounded corners=2mm] plot coordinates {(-0.5, 0.5) (5.6, 0.8) (3.4, -1)};
	\fill [pattern color=blue!80, pattern=north west lines, draw=none, rounded corners=2mm] plot coordinates {(3.4, -1) (-2, -1) (-0.5, 0.5)};
	\fill [black!10, draw=blue, line width=0.5pt, rounded corners=4mm] plot coordinates {(-0.7, 0.29) (5.6, 0.8) (3.2, -1.01)};
	\fill [black!10, draw=black, rounded corners=2mm] plot coordinates {(3.2, -1) (-2, -1) (-0.7, 0.3)};
	\end{scope}
	
	\begin{scope}[shift={(4, 1)}]
	\fill [pattern color=blue!80, pattern=north west lines, draw=blue, line width=0.5pt, rounded corners=2mm] plot coordinates {(0.6, 1) (-1.6, -0.8) (4.5, -0.5)};
	\fill [pattern color=blue!80, pattern=north west lines, draw=none, rounded corners=2mm] plot coordinates {(4.5, -0.5) (6, 1) (0.6, 1)};
	\fill [black!30, draw=blue, line width=0.5pt, rounded corners=4mm] plot coordinates {(0.8, 1.01) (-1.6, -0.8) (4.7, -0.29)};
	\fill [black!30, draw=black, rounded corners=2mm] plot coordinates {(4.7, -0.3) (6, 1) (0.8, 1)};
	\node at (4.3, -2.7) {$(M, g, k)$};
	\end{scope}
	
	\node at (2, 0) {$\bullet$};
	\node at (2.3, -0.3) {$O$};
	
	\draw [black, rounded corners=2mm] plot coordinates {(-4.5, -0.5) (-2, 2) (4.6, 2)};
	\draw [black, rounded corners=2mm] plot coordinates {(8.5, 0.5) (6, -2) (-0.6, -2)};
	
	\begin{scope}[shift={(-4,-1)}]
	\fill [black!10, draw=black] plot [smooth] coordinates {(0, 0) (1, 0.5) (1.3, 1.5) (1.1, 2.5) (2, 2.5) (2.9, 2.5) (2.7, 1.5) (3, 0.5) (4, 0)};
	\end{scope}
	
	\begin{scope}[shift={(4,1)}]
	\fill [black!30, draw=black] plot [smooth] coordinates {(0, 0) (1, 0.5) (1.25, 1.5) (1.1, 2.5) (2, 2.9) (2.9, 2.5) (2.75, 1.5) (3, 0.5) (4, 0)};
	\end{scope}
	
	\end{scope}
	
	\end{scope} 
	\end{tikzpicture}

\end{center}
	\caption{Theorem~A allows to merge an assigned collection of data into an exotic $N-$body solution of the Einstein constraint equations. From~\cite{CarlottoSchoen}, with kind permission of the authors.}
	\label{F16X16.1}
\end{figure}

The construction can be iterated to produce many-body initial data sets.

An alternative gluing construction with \emph{bounded} {sets} $\Omega_a\subset \hyp_a$ has been previously carried-out in~\cite{CCI,CCI2}.

\section{Gluing methods}
 \label{ss8X16.1}

One of the central problems in mathematical general relativity is the construction of solutions to the constraint equations. In spite of an impressive body of work on this by many researchers$^{\mathrm{\eq{Fn12X16.1}}}$
we are still very far from understanding the problem at hand.
Theorem~A  and its technique of proof are  important contributions to the subject, the full consequences of which remain to be explored. The theorem belongs to the category of ``gluing theorems'' for initial data sets, initiated by Justin Corvino in his thesis supervised by Rick Schoen~\cite{Corvino};
compare~\cite{CorvinoSchoen2}. In those last two papers gluing theorems across annuli are developed. The basic ideas behind  Theorem~A  are essentially identical to those of~\cite{Corvino,CorvinoSchoen2} with one key difference: the gluing region is now allowed to be non-compact. This introduces the need to develop new function spaces and establish some key inequalities in those spaces. We will return to those issues in Section~\ref{s8X16.11} below, but before doing this let us shortly discuss some basic aspects of ``gluing''.

Consider, then, two initial data sets 
which are close to each other on a domain $\Omega\subset \hyp$. One further assumes that $\Omega$ has exactly two boundary components, with each component of $\Omega$ separating $\hyp$ into two. The reader can think of $\Omega$ as an annulus in $\hyp=\R^n$, or the region between the cones of Figure~\ref{F11IX16.1} in $\R^n$. The basic idea is to use the inverse function theorem to construct a new initial data set which will coincide with the first initial data set near a component of the boundary, and with the second initial data set near the other component of $\partial \Omega$.

Let us denote by $P$ the linearisation of the map, say $\mathcal C$, which to a pair $(g,K)$ assigns the right-hand sides of \eq{11IX16.3}-\eq{11IX16.4}:
\begin{equation}
 \label{defJrho}
 {\mathcal C}(g,K) :=
\left(
\begin{array}{l}
R(g)-|K|^2 + (\mathrm{tr}  K)^2-2\Lambda\\
2(-\nabla^jK_{ij}+\nabla_i\;\mathrm{tr}  K)
  \\
\end{array}
\right)
\,.
\end{equation}
Let $P^*$ denote the formal adjoint of $P$. A somewhat lengthy calculation gives:
\bea
 \label{4}
  \lefteqn{
  \phantom{xxx}
   P^*(N,Y)
   =
   }
   &&
\\
 \nonumber
  &&
   \left(
\begin{array}{l}
\nabla^lY_l K_{ij}-2K^l{}_{(i}\nabla_{j)}Y_l+
K^q{}_l\nabla_qY^lg_{ij}-\Delta N g_{ij}+\nabla_i\nabla_j N\\
\phantom{xxx}\; +(\nabla^{p}K_{lp}g_{ij}-\nabla_lK_{ij})Y^l-N \mathrm{Ric}(g)_{ij}
+2NK^l{}_iK_{jl}-2N \mathrm{tr}K K_{ij}
\\
 \phantom{xx}
\\
2(\nabla_{(i}Y_{j)}-\nabla^lY_l g_{ij}-K_{ij}N+\mathrm{tr}K\; N g_{ij})
\end{array}
\right)
 \,,
\eea
where the first two lines of the right-hand side should be understood as a single one.

Consider the linearised equation:
\bel{16X16.4}
 P (\delta g, \delta K) = (\delta \rho, \delta J)
 \,.
\ee
Whenever $P P^*$ is an isomorphism, solutions of \eq{16X16.4}
can be obtained by solving the equation
\bel{9X16.11}
 P P^* (N,Y) = (\delta \rho, \delta J)
\ee
for a function $N$ and a vector field $Y$, and setting
$$
  (\delta g, \delta K) = P^*( N, Y)
  \,.
$$
It turns out that $P P^*$ is elliptic in the sense of Agmon, Douglis and Nirenberg~\cite{Morrey}. So the only {essential} obstruction to solving \eq{9X16.11} is the kernel of $PP^*$. Now, because of the form of the operator, the kernel of $PP^*$ {typically} coincides with the kernel of $P^*$, {with the latter contained in Ker$\,PP^*$ in any case}. Nontrivial elements of Ker$\, P^*$ are called \emph{Killing Initial Data}, abbreviated as KIDs. The terminology is due to the fact that KIDs are in one-to-one correspondence with Killing vectors in the vacuum spacetime obtained by evolving the initial data set~\cite{Moncrief75}.
KIDs $(N,Y)$ are thus solutions of the
the \emph{KIDs equation}:
\bea
 \label{2X15.4}
 \nabla_{(i}Y_{j)}
   &= & N K_{ij}
 \,,
\eea
\bea
\lefteqn{\hspace{.2cm}
 \nabla_i \nabla_j N
  =
  }
 &&
\\
\nonumber
 &&
   \bigg(
   \mathrm{Ric}(g)_{ij}
    - 2K^l{}_iK_{jl} +\mathrm{tr}K K_{ij}
     + \frac{1}{ 1-n }
   \big(
    R
    +(\mathrm{tr}K )^2
    -  K^{ql}K_{ql} \big) g_{ij}
 \bigg)
        N
\\
 \nonumber&&
  +\big(\nabla_lK_{ij}+ \frac{1}{n-1}(\nabla^{p}K_{lp}- \nabla_l \mathrm{tr}K)g_{ij}
   \big)Y^l
   + 2K^l{}_{(i}\nabla_{j)}Y_l
 \,,
\eea
In the time-symmetric case with $K\equiv Y\equiv 0$, KIDs are reduced to a single function $N$, and
the KID equations reduce to the \emph{static KIDs equation},
\bea
 \nabla_i \nabla_j N
  =
   \big(
   \mathrm{Ric}(g)_{ij}
     + \frac{R}{ 1-n }
 \big)
        N
 \,.
\eea
%

\subsection{A toy model: divergenceless vector fields}
 \label{ss16X16.13}

To illustrate how this works in a simpler setting, consider the Maxwell constraint equation for a source-free {(that is to say, divergence-free)} electric field $ E$,
\bel{9X16.4}
 P(E):= \mathrm{div} E = 0
 \,.
\ee
The formal adjoint of the divergence operator is the negative of the gradient, so that the ``KID equation'' in this case reads
\bel{9X16.5}
 P^* (u)  \equiv - \nabla u  = 0
 \,.
\ee
The gradient operator has no kernel on a domain with smooth boundary if $u $ is required to vanish on $\partial \Omega$; in fact, the vanishing at a single point of the boundary would suffice. So the equation
\bel{9X16.6}
 \mathrm{div} E = \rho
\ee
can be solved by solving the Laplace equation for $u $,
\bel{9X16.7}
 PP^* (u)  \equiv  - \mathrm{div} \nabla u  \equiv - \Delta u  = \rho
 \,,
\ee
with zero Dirichlet data.

Consider, then, the following toy problem:

\begin{Problem}
  \label{Probl9X16.1}
Let $E_i$, $i=1,2$, be two source-free electric fields on $\R^n$. Find a source-free electric field $E$ which coincides with $E_1$ on a ball $B(R_1)$ of radius $R_1$ and coincides with $E_2$ outside a ball of radius $R_2>R_1$.
\end{Problem}

Note that if $E_2\equiv 0$, and if we can solve Problem~\ref{Probl9X16.1}, we will have screened away the electric field $E_1$ without introducing any charges in the system: this is the screening of the electric field with an electric field. We will also have constructed an infinite dimensional space of compactly supported divergence free vector fields as $E_1$ varies, with complete control of $E$ in $B(R_1)$.

Note that the alternative construction of the space of such vector fields, by setting $E=\mathrm{curl}\, A$ for some compactly supported vector field $A$, does not provide explicit control of $E|_{B(R_1)}$.

As a first attempt to solve the problem, let $\chi$ be a radial cut-off function which equals one near the sphere $S(R_1)$ and which equals zero near $S(R_2)$. Set
$$
 E_\chi = \chi E_1 + (1-\chi)E_2
 \,.
$$
Since both $E_i$ are divergence-free we have
$$
 \rho_\chi:= \mathrm{div} E_\chi = \nabla \chi \cdot (E_1-E_2)
  \,,
$$
and there is no reason for $  \rho_\chi$ to vanish. However, if $u $ solves the equation
\bel{9X16.71}
 PP^* (u)  \equiv  - \Delta u  = - \rho_\chi
\ee
with vanishing boundary data, then
\bel{9X16.72}
 E = E_\chi + P^*(u )
\ee
will be divergence free:
$$
 \mathrm{div} E = \mathrm{div} ( E_\chi +  P^*(u ))
  = \rho_\chi + P P^*(u ) =0
  \,.
$$
Now, on $S(R_i)$ we have
\bel{9X16.15}
 E|_{S(R_i)} = E_i - \nabla u
 \,,
\ee
and there is no reason why this should coincide with $E_i$. We conclude that this approach fails to solve the problem.

{Replacing Dirichlet data by Neumann data will only help if both $E_i|_{S(R_i)}$ are purely radial, as suitable Neumann data will only guarantee continuity of the normal components of $E$.}

It turns out that there is  trick to make this work in whole {generality:   modify \eq{9X16.7} by introducing weight-functions $\psi$ which vanish very fast at the boundary.} An example, which will lead to solutions on the annulus which can be extended to the whole of $\R^n$ in a high-but-finite differentiability class, is provided by the functions
\bel{9X16.16}
 \psi=(r-R_1)^\sigma(R_2-r)^\sigma
 \,,
  \quad r\in (R_1,R_2)
  \,,
\ee
with  some large positive number $\sigma$. Another useful example, which will lead to smoothly-extendable solutions, is
\bel{9X16.17}
 \psi=(r-R_1)^\alpha(r-R_2)^\alpha \exp\Big(-\frac s {(r-R_1)(R_2-r)}\Big)
 \,,
  \quad r\in (R_1,R_2)
 \,,
\ee
with $\alpha \in \R$ and $s> 0$. 
(The prefactors involving $\alpha$ in \eq{9X16.17} are useful when constructing a consistent functional-analytic set-up, but are essentially irrelevant as far as the blow-up rate of $\psi$ near the $S(R_i)$'s is concerned.)

As such,  instead of \eq{9X16.7} consider the equation
\bel{9X16.73}
 P(\psi^{2}P^* (u) ) \equiv  - \mathrm{div} (\psi^{2} \nabla  u)  = - \rho_\chi
 \,.
\ee
Solutions of \eq{9X16.73} could provide a solution of Problem~\ref{Probl9X16.1} if one replaces \eq{9X16.72} by
\bel{9X16.+73}
 E = E_\chi + \psi^{2} P^*(u )= E_\chi - \psi^{2} \nabla u
 \,.
\ee
Now, solutions of \eq{9X16.73} are, at least formally, minima of the functional
\bel{9X16.51}
 I= \int_\Omega \frac 12 \psi^2 |\nabla u|^2 +\rho_\chi u
 \,.
\ee
Supposing that minimisation would work, one will then obtain a solution $u$ so that $\psi \nabla u$ is in $L^2$. Since $\psi$ goes to zero at the boundary very fast, $\nabla u$ is likely to blow up. In an ideal  world, in which ``$L^2$'' is the same as ``bounded'', $\nabla u$ will behave as $\psi^{-1}$ near the boundary. The miracle is that \eq{9X16.73} involves $\psi^2 \nabla u$, with one power of $\psi$ spare, and so the derivatives of $u$ would indeed tend to zero as $\partial \Omega$ is approached.

This naive analysis of the boundary behaviour turns out to be essentially correct: choosing the exponential weights \eq{9X16.17}, $\psi^2 \nabla u $ will extend smoothly by zero across the boundaries when the $E_i$'s are smooth. A choice of power-law weights \eq{9X16.16} will lead to extensions of differentiability class determined by the exponent $\sigma$, with a loss of a finite number of derivatives due to the fact that $L^2$ functions are not necessarily bounded, and that there is a loss of differentiability when passing from Sobolev-differentiability to classical derivatives.

It then remains to show that minimisation works in a carefully chosen space. This requires so-called ``coercitivity inequalities''. For the functional \eq{9X16.51} the relevant inequality is the following \emph{weighted Poincar\'e inequality:}
\bel{9X16.52}
 \int_\Omega \psi^2 |u|^2 \le C \int_\Omega \varphi^2 \psi^2 |\nabla u|^2
 \,.
\ee
Here $\varphi = (r-R_1)^{-1}(R_2-r)^{-1}$ when  $\psi$ is given by \eq{9X16.16} with $\sigma \ne -1/2$, and
$\varphi =(r-R_1)^{-2}(R_2-r)^{-2}$ for the exponential weights  $\psi$ given by \eq{9X16.17} with $s\ne 0$, cf.\ e.g.\ \cite{ChDelay,Corvino}.

There is an obvious catch here, namely \eq{9X16.52} cannot possibly be true since it is violated by constants. However, \eq{9X16.52} holds on the subspace, say $F$, of functions which are $L^2$-orthogonal to constants, when using the weights described above. This turns out to be good enough for solving Problem~\ref{Probl9X16.1}. For then one can carry out the minimisation on $F$, finding a minimum $u\in F$. The function $u$ will solve the equation up to an $L^2$-projection of the equation on constants. In other words, we will have
\bel{9X16.41}
 \int_\Omega f ( - \mathrm{div} (\psi^2 \nabla u) + \rho_\chi) = 0
 \,,
\ee
for all differentiable $f$ such that $\int_\Omega f = 0$.
Now, integrating the equation \eq{9X16.73} against the constant function $f\equiv 1$ we find
\beal{10X16.3}
 \lefteqn{
 \int_\Omega  (  \mathrm{div} (\psi^2 \nabla u) - \rho_\chi)
   =
 \int_\Omega    \mathrm{div} (\psi^2 \nabla u - E_\chi)
 }
 &&
\\
 \nonumber
 &  &=
 \int_{\partial \Omega} (\psi^2 \nabla u - E_\chi) \cdot m = \int_{S(R_1)} E_1 \cdot n- \int_{S(R_2)} E_2 \cdot n
 \,,
\eea
where $m$ is the outer-directed normal to $\partial \Omega$, and $n$ is the radial vector $\vec x/|\vec x|$. Since the $E_i$'s are divergence-free it holds that
$$
 \int_{S(R_i)} E_i \cdot n = \int_{B(R_i)} \mathrm{div}(E_i)  =0
 \,.
$$
It follows that the right-hand side of \eq{10X16.3} vanishes, hence \eq{9X16.41} holds for all differentiable functions $f$, and  $u$ is in fact a solution of \eq{9X16.73}.

(Strictly speaking, when solving \eq{9X16.73} by minimisation using the inequality \eq{9X16.52}, Equation~\eq{9X16.73} should be replaced by
\bel{9X16.7+3}
  \mathrm{div} ( \varphi^2\psi^{2} \nabla  u)  =  \rho_\chi
 \,.
\ee
This does not affect the discussion so far, and only leads to a shift of the powers $\alpha$ and $\sigma$ in \eq{9X16.16}-\eq{9X16.17}.)

We conclude that the answer to Problem~\ref{Probl9X16.1} is yes, as already observed in~\cite{CorvinoPollack,ErwannTT}.

One can likewise solve variants of Problem~\ref{Probl9X16.1} with gluing regions which are not annuli, e.g. a difference of two coaxial cones with distinct apertures as in the Carlotto-Schoen Theorem~A. If one of the glued vector fields is taken to be trivial, one obtains configurations where the electric field extends all the way to infinity in open cones and vanishes in, e.g., a half-space.

The gluing construction for the linearised relativistic constraint equations proceeds essentially in the same way. There, in addition to the weighted Poincar\'e inequality \eq{9X16.52} one also needs a \emph{weighted Korn inequality} for vector fields $X$:
\bel{9X16.52+}
 \int_\Omega \psi^2 |X|^2 \le C \int_\Omega \varphi^2 \psi^2 |S(X)|^2
 \,,
\ee
where  $S(X)$ is the symmetric two-covariant vector field defined as
\bel{16X16.5}
 S(X)_{ij} = \frac 12 (\nabla_i X_j + \nabla_j X_i)
 \,.
\ee
In \eq{9X16.52+} one needs to assume $\sigma \not\in\{ -n/2,-n2/-1\}$ for the weight \eq{9X16.16}, and $s\ne 0$ when $\psi$ is given by \eq{9X16.17}. If the metric $g$ has non-trivial Killing vectors, which are solutions of the equation $S(X)=0$, then  \eq{9X16.52+} will hold for vector fields $X$ which are in a closed subspace transverse to  the space of Killing vectors, with a constant depending upon the subspace.

Because the space of KIDs is trivial for generic metrics~\cite{BCS}, the problem of solving modulo kernel does not arise in generic situations.

The full non-linear gluing problem for the scalar curvature, or for vacuum initial data, is solved using the above analysis of the linearised equation together with a tailor-made version of the inverse function theorem. The reader is referred to \cite{ChDelay,Corvino} for details. There it is also shown how to treat problems where existence of KIDs cannot be ignored.

\subsection{Previous generalisations}
 \label{ss16X16.12}

As already mentioned, the original papers~\cite{Corvino,CorvinoSchoen2} were concerned with gluing across an annulus in $\R^n$, with vanishing cosmological constant. This has been generalised to various other set-ups, and used to prove the following:

\begin{enumerate}
 \item {In~\cite{LiYu} the method is used to construct non-trivial black hole spacetimes with smooth asymptotic structure.}
      \item In~\cite{ChDelay2,CorvinoAHP} the gluing method was used to construct  large families of vacuum spacetimes  which are asymptotically simple in the sense of Penrose~\cite{penrose:asymptotic}.
  \item In~\cite{ChDelay} it has been shown how to reduce the gluing problem to the verification of a few properties of the weight functions $\varphi$ and $\psi$, and to the verification of the Poincar\'e and Korn inequalities. It has also been shown how to use the technique to control the asymptotics of  solutions in asymptotically flat regions.
      \item In~\cite{CIP:CMP} the method has been used to carry out localised vacuum connected-sums of vacuum initial data sets.
      \item In~\cite{ChDelayAH} the method is used to construct constant negative scalar curvature metrics with exact Schwarzschild-anti de Sitter asymptotics.
      \item It has been shown in~\cite{CPP,ChPollack,CortierKdS} how to use the gluing method to construct manifolds with constant positive scalar curvature containing periodic asymptotic ends.
          \item The differentiability thresholds for the applicability of the method have been lowered in~\cite{ChDelayHilbert}. It has also been shown there how to construct a Banach manifold structure for the set of vacuum initial data under various asymptotic conditions using the general ideas developed in the process of gluing.
                  \item In~\cite{ErwannInterpolating} gluings are done by interpolating scalar curvature.
              \item In~\cite{ErwannTT} the gluing method has been used to construct compactly supported solutions for a wide class of underdetermined elliptic systems. {As a particular case, for any open set $\mycal U$ one obtains an infinite-dimensional space of solutions of the vector constraint equation which are compactly supported in $\mycal U$.}
                  \item In~\cite{DelayMazzieri} the gluing is used to make local ``generalised connected-sum'' gluings along submanifolds.
                      \item As already pointed-out, in~\cite{CCI,CCI2} the gluing is used to construct initial data for the many-body problem in general relativity.
                  \item In~\cite{CorvinoHuang} the gluing method is used to deform initial data satisfying the dominant energy condition to ones where the condition is strict.
                          \item In~\cite{ChBieri} the gluing is used {to remove the previous smallness conditions in the assertion that maximal globally hyperbolic developments of asymptotically flat initial data sets contain null hypersurfaces with generators complete to the future.}
\end{enumerate}

\section{Further elements of the proof}
 \label{s8X16.11}

As already mentioned, the main idea of the proof of Theorem~A  is essentially identical to that of the Corvino-Schoen gluing, summarised above. There is, however, a significant amount of new work involved, {which fully justifies the length of the paper and  publication in \emph{Inventiones}.}

The first extraordinary insight is to imagine that the result can be true at all when $\Omega$ is the difference of two cones with different apertures, smoothed out at the vertex, see Figure~\ref{F9X16.1}.
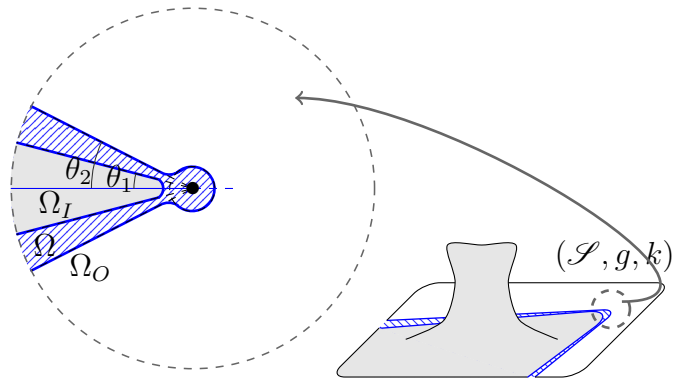
\begin{figure}[ht!]
\begin{center}
\begin{tikzpicture}
	
	\begin{scope}[scale=0.5, shift={(20, -4)}]
	\node at (5.4, 1) (MagSmall) {};
	\fill [pattern color=blue!80, pattern=north west lines, draw=blue, line width=0.5pt, rounded corners=2mm] plot coordinates {(-0.5, 0.5) (5.6, 0.8) (3.4, -1)};
	\fill [pattern color=blue!80, pattern=north west lines, draw=none, rounded corners=2mm] plot coordinates {(3.4, -1) (-2, -1) (-0.5, 0.5)};
	\fill [black!10, draw=blue, line width=0.5pt, rounded corners=4mm] plot coordinates {(-0.7, 0.29) (5.6, 0.8) (3.2, -1.01)}; 
	\fill [black!10, draw=black, rounded corners=2mm] plot coordinates {(3.2, -1) (-2, -1) (-0.7, 0.3)};         
	\draw [black, rounded corners=2mm] plot coordinates {(-0.5, 0.5) (0.5, 1.5) (7, 1.5) (4.5, -1) (3.4, -1)};          
	\fill [black!10, draw=black] plot [smooth] coordinates {(0, 0) (1, 0.5) (1.3, 1.5) (1.1, 2.5) (2, 2.5) (2.9, 2.5) (2.7, 1.5) (3, 0.5) (4, 0)};
	\draw [black!60, line width=1pt, dashed] (5.4, 0.8) circle (0.5);
	\node at (5.5, 2.2) {${\small (\hyp, {g}, {k})}$};              
	\end{scope}
	
	\begin{scope}[scale=0.3, shift={(24, 0)}]
	
	\clip (0, 0) circle (8);
	\node at (4, 4) (MagLarge) {};
	
	\begin{scope}
	\filldraw [pattern color=blue!80, pattern=north east lines, draw=blue, line width=1pt] (-8, 4.1) -- (-1.3, 0.66) arc [radius=0.75, start angle=-110, delta angle=50] -- (-0.69, 0.71) arc [radius=0.98, start angle=132, delta angle=-264] -- (-0.69, -0.73) arc [radius=0.75, start angle=60, delta angle=60] -- (-8, -4.1);
	\filldraw [fill=black!10, draw=blue, line width=1pt] (-8.02, 2.105) -- (-1.5, 0.4) arc [radius=0.5, start angle=53, delta angle=-106] -- (-1.5, -0.4) -- (-8.02, -2.105); 
	\end{scope}
	
	\begin{scope}
	\draw [black, dashed] plot coordinates {(-8, 4.1) (0, 0) (-8, -4.1)};
	\draw [black, dashed] plot coordinates {(-8, 2.1) (0, 0) (-8, -2.1)};
	\draw [black, dashed] (-0.03, -0.02) circle (0.99);
	\end{scope}
	
	\begin{scope}
	\draw [blue, ->] plot coordinates {(-8, 0) (-0.1, 0)};
	\draw [blue, dashed] plot coordinates {(0, 0) (2, 0)};
	\node at (0, 0) {$\bullet$};
	\node at (0.35, 0.35) {$~$};
	\end{scope}
	
	\begin{scope}
	\draw [black!80] (-4, 2.05) arc [radius=4.8, start angle=155, delta angle=25];
	\node at (-5, 0.8) {${ \theta_2}$};
	\draw [black!80] (-2.5, 0.656) arc [radius=2.59, start angle=165, delta angle=15];
	\node at (-3.2, 0.45) {$\theta_1$};
	\end{scope}
	
	\begin{scope}
	\node at (-6, -0.7) {$\Omega_I$};
	\node at (-6.5, -2.5) {$\Omega$};
	\node at (-4.5, -3.5) {$\Omega_O$};
	\end{scope}
	
	\draw [black!60, line width=1pt, dashed] (0, 0) circle (7.99);
	
	\end{scope}
	\path[every node/.style={font=\sffamily\small}, ->, black!60, line width=1pt] (MagSmall) edge [out=0, in=0] (MagLarge);
	
\end{tikzpicture}
\end{center}
	\caption{\label{F9X16.1}Regularized cones and the gluing region $\Omega$, from \cite{CarlottoSchoen} with kind permission of the authors. $\Omega_I\subset \Omega$ is the interior cone, $\Omega_O$ is the region outside the larger cone.}
\end{figure}
This is the geometry that we are going to assume in the remainder of this section.

Next, all generalisations of~\cite{Corvino,CorvinoSchoen2} listed in Section~\ref{ss16X16.12} involve gluing across a compact boundary. In the  current case $\partial \Omega$ is not compact, and so some analytical aspects have to be revisited.  In addition to weights governing decay at $\partial \Omega$, radial weights need to be introduced in order to account for the infinite extent of the cone.

Let $d(p)$ denote the distance from $p\in \Omega$ to $\partial \Omega$. Let  $\theta_1<\theta_2$ be the respective apertures of the inner and outer cones and let $\theta$ denote the angle away from the axis of the cones. A  substantial part of the paper consists in establishing inequalities in the spirit of \eq{9X16.52} and \eq{9X16.52+}  with $\varphi=r$ and with weight functions $\psi$ which are smooth everywhere, behave like $d^\sigma$ for small $d$, and are equal to
\bel{9X16.81}
 \psi = r^{n/2-q} (\theta-\theta_1)^\sigma(\theta_2-\theta)^\sigma
 \,,
\ee
for large distances, with $q,\sigma>0$. More precisely,
let $\sigma>0$ be large enough, and assume that $0<q< (n-2)/2$, with $q\ne (n-4)/2$ for $n\ge 5$. Suppose that $g$ is the Euclidean metric and let $\Omega$ be as above. Let $ \phi$ be  a positive function which for large distances equals $\theta-\theta_1$ and $\theta_2-\theta$ close to the inner and outer cones, respectively,
 and which behaves as the distance from $
\partial \Omega$ otherwise. Then there exists a constant $C$ such that for all differentiable functions $u$ and vector fields $X$, both with bounded support in $\Omega$ (no conditions at $\partial\Omega$),
the following inequalities are true:
\beal{11Xi16.1}
 &&\int_\Omega |u|^2 r^{-n+2q}  \phi^\sigma \le C
 \int_\Omega |\nabla u |^2 r^{2-n+2q}  \phi^\sigma
 \,,
%
%
\\
 &&\int_\Omega |Y|^2 r^{-n+2q}  \phi^\sigma \le C
 \int_\Omega |S(Y)|^2 r^{2-n+2q} \phi^\sigma
 \,.
\eeal{11Xi16.4}

A clever lemma relying on the coarea formula, (\cite[Lemma~4.1]{CarlottoSchoen}), reduces the proof of the inequalities \eq{11Xi16.1}-\eq{11Xi16.4}
to the case  $\phi\equiv 1$.

A key point in the proof of \eq{11Xi16.4}  is the inequality  established in  \cite[Proposition~4.5]{CarlottoSchoen} (also known in  \cite{CarlottoSchoen} as ``Basic Estimate II''), which takes the form
\bel{12X16.3}
\int_{\Omega}|\nabla Y|^2 r^{2-n+2q}  \phi^\sigma \leq C\int_{\Omega}|S(Y)|^2r^{2-n+2q}  \phi^\sigma
 \,.
\ee
The justification of \eq{12X16.3} requires considerable ingenuity.

It is simple to show that \eq{11Xi16.1}-\eq{11Xi16.4} continue to hold for asymptotically Euclidean metrics which are close enough to the Euclidean one, {with uniform constants}. As explained in Section~\ref{ss16X16.13}, these inequalities provide the stepping stones for the analysis of the linear equations.

Note that the radial weights in \eq{11Xi16.1}
guarantee that affine functions are not in the space obtained by completing $C_c^1(\overline \Omega)$ with respect to the {norm defined by the  right-hand side}. A similar remark concerning vectors with components which are affine functions applies in the context of \eq{11Xi16.4}.
This guarantees that neither  KIDs, nor asymptotic KIDs, interfere with the construction, which would otherwise have introduced a serious obstruction to the argument.

Once these \emph{decoupled} functional inequalities are gained, a perturbation argument ensures coercivity of the adjoint linearised constraint operator (in suitable doubly-weighted Sobolev spaces). This allows one to use direct methods to obtain existence of a unique global minimum for the functional whose Euler-Lagrange equations are  the linearized constraints. We refer the reader to Propositions 4.6 and  4.7 of~\cite{CarlottoSchoen} for precise statements.

The argument of~\cite{CarlottoSchoen} continues with a Picard iteration scheme, which allows one to use the analysis of the linear operator to obtain solutions to the nonlinear problem under a smallness condition.  {This is not an off-the-shelf argument: it involves some  delicate  choices of functional spaces for the iteration, where one takes a  combination  of weighted-Sobolev and weighted-Schauder norms. Alternatively one could use~\cite[Appendix~G]{ChDelay2} at this stage of the proof, after establishing somewhat different estimates, compare~\cite{ChDelayExoticAERiemannian}.}

To end the proof it suffices to start moving the cones to larger and larger distances in the asymptotic region, so that the metric on $\Omega$ approaches the flat one.  When the tips of the cones are far enough the smallness conditions needed to make the whole machinery work are met, and Theorem~A  follows.

An interesting, and somehow surprising, aspect of the result is the fact that \emph{no matter how small the cone angles are}, the ADM energy-momentum of the glued data provides an arbitrarily good approximation of the ADM energy-momentum of the given data when the vertex of the cones is chosen far enough in the asymptotic region. This is proven in Section 5.6 of~\cite{CarlottoSchoen} and is then exploited in the construction of $N$-body Carlotto-Schoen solutions, already presented in Section~\ref{ss16X16.11}. This is the object of Section 6 of their paper.

\section{Beyond Theorem~A}
 \label{s8X16.22}

The results of Carlotto and Schoen have meanwhile been extended in a few directions.

In the initial-data context, gluings in the same spirit have been done in~\cite{ChDelayExotic} for asymptotically hyperbolic initial data sets. In terms of the half-space model for hyperbolic space,  the analogues of cones are
half-annuli extending to the conformal boundary at infinity. As a  result one obtains e.g.\ non-trivial constant scalar curvature metrics which are exactly hyperbolic in half-balls centered at the conformal boundary.
{We provide more details in Section~\ref{ss11X16.5} below.}

In a Riemannian asymptotically Euclidean setting, with $K_{ij}\equiv 0$ so that only the scalar curvature matters, the {following  generalisations are straightforward}:

\begin{enumerate}
 \item Rather than gluing an asymptotically Euclidean metric to a flat one, any two asymptotically metrics $g_1$ and $g_2$ are glued together.
 \item In the spirit of~\cite{ErwannInterpolating}, the gluings at zero-scalar curvature can be replaced by gluings where the scalar curvature of the final metric equals
     $$
      \chi R(g_1) + (1-\chi) R(g_2)
     $$
     where, as before, $\chi$ is a cut-off function varying between zero and one in the gluing region. Thus, the scalar curvature of the final metric is sandwiched between the scalar curvatures of the original ones. This reduces of course to a zero-scalar-curvature gluing if both $g_1$ and $g_2$ are scalar-flat.
      \item The geometry of the gluing region can be allowed to be somewhat more general than the interface between two cones~\cite{ChDelayExoticAERiemannian}.
\end{enumerate}

{A few more details about this can be found in Section~\ref{ss11X16.6} below.}

\subsection{Asymptotically hyperbolic gluings}
 \label{ss11X16.5}

Let us describe here one of the gluing constructions in~\cite{ChDelayExotic}, the reader is referred to that reference for some more general ``exotic hyperbolic gluings''. The underlying manifold   is taken  to be the ``half-space model'' of hyperbolic space:
$$
 \mcH
  =\{(z,\theta)|\ z>  0, \theta \in\R^{n-1}\}\subset \R^n
  \,.
$$
One wishes to  glue together metrics asymptotic to each other while interpolating their respective
scalar curvatures.
The first metric is assumed to take the form,  in suitable local coordinates,
\bel{bgenerale}
 g =\frac {1} {z^2}
  \big(
   {(1+O(z))dz^2 + \underbrace{h_{AB}(z,\theta^C)d\theta^A d\theta^B}_{=:h(z)}}
    +O(z)_Adz\, d\theta^A
   \big)
 \,,
\ee
where $h(z)$ is a continuous family of Riemannian metrics on $\R^{n-1}$.

We define
$$
 B_\lambda  := \{z>0\,,\    \underbrace{\sum_i(\theta^i)^2}_{=:|\theta|^2} + z^2 <\lambda^2\}
 \,,\;\;\;A_{\epsilon,\lambda}=B_{\lambda}\setminus \overline{B_\epsilon}
 \,.
$$
The gluing construction will  take place in the region
\bel{29II16.1}
 \Omega=A_{1,4}
 \,.
\ee

Let $\hg$ be a second metric   on $B_5$ which is close to $ g$ in $C^{k+4}_{1,z^{-\sigma}}(A_{1,4})$. Here,
for $\phi$ and $\varphi$  --- smooth strictly positive functions
on M, and for $k\in\N$,
we define
$C^{k }_{\phi,\varphi}$ to be the space of $C^{k }$
functions or tensor fields  for which the norm
$$
\begin{array}{l}
\|u\|_{C^{k}_{\phi,\varphi}(g)}=\sup_{x\in
M}\sum_{i=0}^k
\|\varphi \phi^i \nabla^{(i)}u(x)\|_g
\end{array}
$$
is finite.

Let $\chi$ be a smooth non-negative function on  $\mcH$, equal to $1$  on $\mcH\setminus B_3$, equal to zero
on $B_2$, and positive on $\mcH\backslash \overline{B_2}$. We set
\bel{21IX15.1}
g_\chi:=\chi \hg+(1-\chi) g
 \,.
\ee
In~\cite{ChDelayExotic}  a gluing-by-interpolation of the constraint equations is carried out. In the time-symmetric case, the main interest is that of constant scalar curvature metrics, which then continue to have constant scalar curvature, or for metrics with positive scalar curvature, which then remains positive.  Since the current problem is related to the construction of initial data sets for Einstein equations, in general-relativistic matter models  such as   Vlasov  or dust, an  interpolation of scalar curvature is  of direct interest.

One has~\cite{ChDelayExotic}:

\begin{theo}
 \label{propRcassimple}
Let $n/2<k<\infty$, $b\in[0,\frac{n+1}2]$, $\sigma>\frac{n-1}2+b$, suppose that $g-\zg \in {C^{k+4}_{1,z^{-1}}}$.
For all 
$\hg$ close enough to $ g$ in $C^{k+4}_{1,z^{-\sigma}}(A_{1,4})$
there exists a
two-covariant symmetric tensor field $h$ in $C^{k+2-\lfloor n/2\rfloor}(\mcH)$, vanishing outside of $A_{1,4}$, such that the tensor field $g_\chi+h$ defines an asymptotically hyperbolic metric satisfying
\bel{solmodker}
 R(g_\chi+h)=\chi R(\hg)+(1-\chi)R( g)
 \,.
\ee
\end{theo}

A similar result is established in~\cite{ChDelayExotic} for the full constraint equations.

The proof involves ``triply weighted Sobolev spaces'' on $A_{1,4}$ with
$$
 \varphi =x/\rho\,, \quad
 \psi =x^a z^b \rho^c
  \,,
$$
where $a$ and $c$ are chosen large as determined by $k$, $n$ and $\sigma$.
Here $z$ is the coordinate of \eq{bgenerale}, the function $x$ is taken to be any smooth function on $ \Omega$ which equals the $z^2\mathring g$-distance
to
$$\{|\theta|^2 + z^2=1\}\cup \{|\theta|^2 + z^2=4 \}
$$
near this last set, while $\rho:= \sqrt{x^2+z^2}$. The heart of the proof is the establishing of the relevant Poincar\'e and Korn inequalities. Once this is done, the scheme of the proof follows closely the arguments described so far.

One would like to have a version of Theorem~\ref{propRcassimple} with weights which exponentially decay as $x$ tends to zero. However, the triply-weighted Korn inequality needed for this has not been established so far.

\subsection{Asymptotically Euclidean scalar curvature gluings by interpolation}
 \label{ss11X16.6}

We finish this \emph{s\'eminaire} by describing a straightforward generalisation  of Theorem~A in the time-symmetric asymptotically Euclidean setting.

Let $S(p,R)\subset \R^n$ denote a sphere of radius $R$ centred at $p$.
Let $\Omega \subset \R^n$ be a domain with smooth boundary (thus $\Omega$ is open and connected, and $\partial \Omega =\overline \Omega \setminus \Omega$ is a smooth manifold).
We further assume that $\partial \Omega$ has \emph{exactly two} connected components, and that
there exists $R_0\ge 1$ such that
\bel{7X16.24}
 \Omega_{S}:= \Omega \cap S(0,R_0)
\ee
also has exactly two connected components, with
\bel{7X16.25}
 \Omega \setminus B(0,R_0)=\{\lambda p \ | \ p \in \Omega_{S}\,, \lambda \ge 1\}
 \,.
\ee
%

The regularised differences of cones in the Carlotto-Schoen  gluings  provide examples of such sets. Another example is displayed in Figure~\ref{F15X16.2}.
\begin{figure}[th]
  \centering
    \includegraphics[scale=.3]{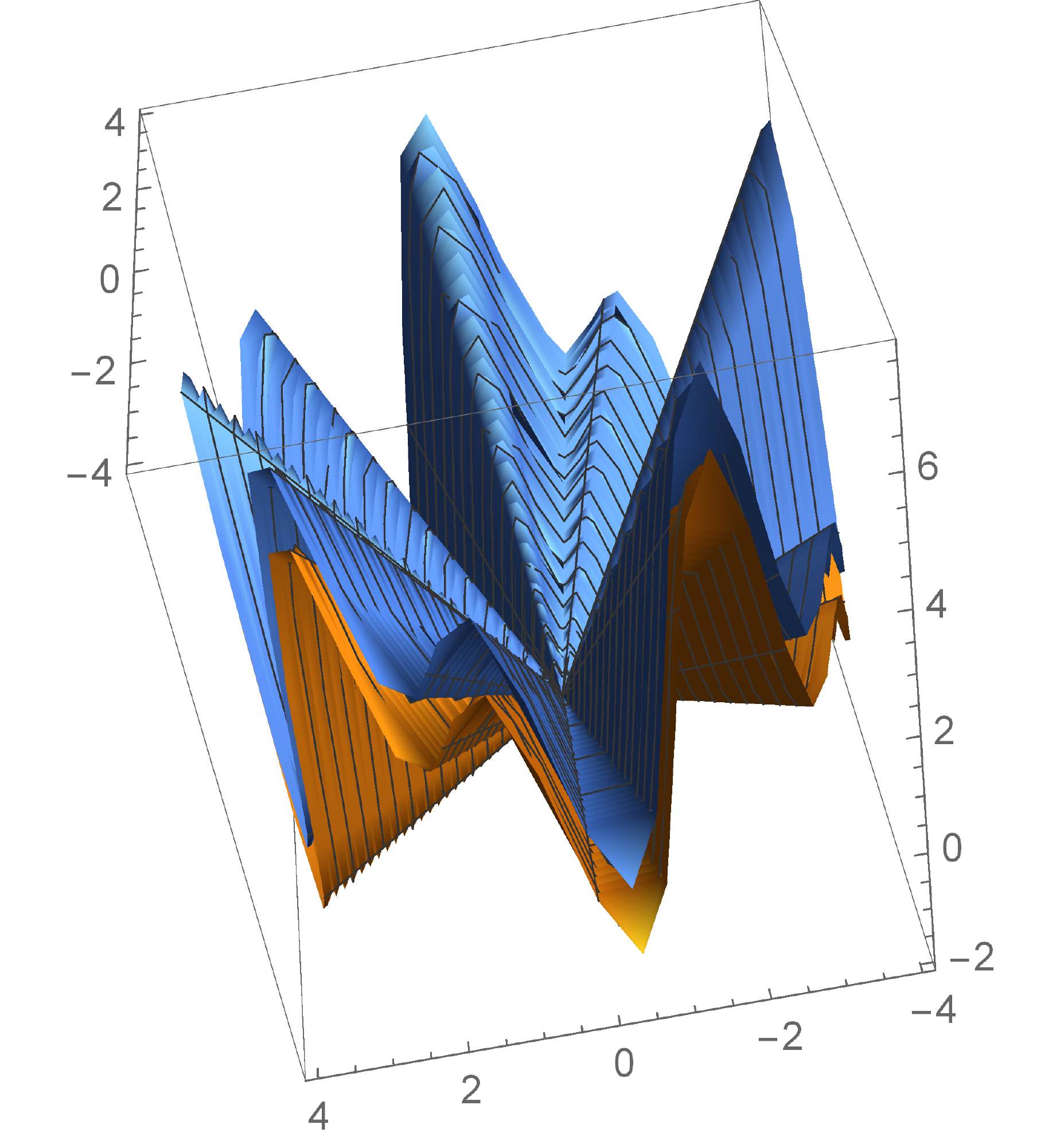}
\caption{A possible set $\Omega$,  located between the two surfaces. The {glued} metric coincides with $g_1$ above the {higher} surface, {coincides} with $g_2$ below the {lower surface}, and is scalar-flat if both $g_1$ and $g_2$ were.
    \label{F15X16.2}}
\end{figure}
Since the construction can be iterated, the requirement that $\Omega$ be connected is irrelevant.

We let $x:\overline \Omega\to \R$ be any smooth defining function for $\partial \Omega$ which has been chosen so that
\bel{14VI14.1}
 \mbox{$x(\lambda p) = \lambda x(p)$ for $\lambda \ge 1$ and for  $p, \lambda p \in \Omega\setminus B(0,R_0)$.}
\ee
Equivalently, for $p\in\Omega_S$ and $\lambda$ larger than one, we require $x(\lambda p) = \lambda x_S(p)$, where $x_S$ is a defining function for  $\partial\Omega_S$ within $S(0,R_0)$.

We will denote by $r$  a  smooth positive function which coincides with $|\vec x|$ for $|\vec x |\ge 1$.

By definition of $\Omega$ there exists a constant $c$ such that the distance function $d(p)$ from a point $p\in\Omega$ to $\partial \Omega$ is smooth for all $d(p)\le c r(p)$. The function $x$ can be chosen to be equal to $d$ in that region.

For $\beta, s, \mu \in \R$ we define
\bel{5X16.11}
 \varphi=\left(\frac{x}{r}\right)^2r=\frac{x^2}r
 \,,
  \quad
  \psi = \funnyr^{-n/2-\beta}\left(\frac{x}{r}\right)^\sigma e^{-sr/x}=:\funnyr^{\mu}{x}^\sigma e^{-sr/x}
\ee
on $\Omega$.  One can show that the weighted Poincar\'e ineqality \eq{9X16.52} holds with these weights, modulo a {supplementary integral of $|u|^2$ on a compact subset of $\Omega$, for all tensor fields $u$ compactly supported in $\Omega$ as long as $s\ne 0$ and $\beta\equiv \sigma+\mu+n/2\ne 0$. The supplementary integral does not lead to any new difficulties in the proof, which proceeds as described above.}

 In order to carry out the scalar-curvature interpolation, recall that $\partial \Omega$ has exactly two connected components. We denote by $\chi$ a smooth function with the following properties:
\begin{enumerate}
 \item $0\le \chi \le 1$;
   \item $\chi$ equals one in a neighborhood of one of the components and equals zero in a neighborhood of the other component;
       \item on $\Omega\setminus B(0,R_0)$ the function $\chi$ is a function of $x/r$ wherever it is not constant.
\end{enumerate}
The metric $g_\chi$ 
is then defined as in \eq{21IX15.1}.


Letting $\Omega$, $x$ and $r$ be as just described, with $\psi$, $\varphi$ given by \eq{5X16.11}, in~\cite{ChDelayExoticAERiemannian} the following is proved:

\begin{theo}
 \label{ThefulltheoremAE}
Let $\epsilon>0$, $k>n/2$, $\beta\in[-(n-2),0)$, and $\tbeta  <\min(\beta,-\epsilon)$. Suppose that $g-\delta \in  C^{k+4}_{r,r^{\epsilon}} \cap C^\infty$, where $\delta$ is the Euclidean metric.
For all real numbers $\sigma$ and $s>0$  and
$$
 \mbox{all smooth metrics
    $\hg $ close enough to $ g$ in $ C^{k+4}_{r,r^{-\tbeta  }}(\Omega )$}
$$
{there exists on $\Omega$ a unique smooth two-covariant symmetric tensor field $\deltag $
%
%
such that the metric $ g_\chi+\deltag  $ solves}
\bel{fullcolle}
 \sourcesR
    \left[g_\chi+\deltag \right]= \chi\sourcesR(\hat g) +(1-\chi)
    \sourcesR(g)
 \,.
\ee
The   tensor field $ \deltag  $ vanishes at $\partial \Omega$ and can be $C^{\infty}$-extended by zero across $\partial \Omega$, leading to a smooth  asymptotically Euclidean metric $g_\chi+\deltag $.
\end{theo}

There is little doubt that there is an equivalent of Theorem~\ref{ThefulltheoremAE} in the full initial-data context. In fact, the only missing element of the proof at the time of writing of this review is a doubly-weighted Korn inequality with weights as in \eq{5X16.11}.

The smallness assumptions needed in the theorem can be realised by moving the set $\Omega$ to large distances, as in the Carlotto-Schoen theorem. When $\Omega$ does not meet $S(0,1)$, an alternative is provided by ``scaling $\Omega$ up'' by a large factor. This is equivalent to keeping $\Omega$ fixed and scaling-down the metrics from large to smaller distances. It should be clear that the metrics will approach each other, as well as the flat metric, when the scaling factor becomes large.

\bigskip

\noindent{\sc Acknowledgements:} I am grateful to Alessandro Carlotto, Justin Corvino and Erwann Delay for useful comments on iterated drafts of this review. Supported in part by the Austrian Research Fund (FWF): Project P 29517-N16.

\bibliographystyle{smfplain}

\bibliography{../references/netbiblio,
../references/newbiblio,
../references/hip_bib,
../references/reffile,
../references/newbiblio4,
../references/newbiblio2,
../references/besse2,
../references/bibl,
../references/howard,
../references/bartnik,
../references/myGR,
../references/newbib,
../references/Energy,
../references/chrusciel,
../references/dp-BAMS,
../references/prop2,
../references/besse
}

\end{document}